\documentclass[10pt]{article}
\usepackage{amssymb}
\usepackage{amscd}
\usepackage{amsthm}
 \usepackage[utf8]{inputenc}
\usepackage{stmaryrd}
\usepackage{hyperref}
\newtheorem{theorem}{Theorem}[section]

\newtheorem{e-proposition}[theorem]{Proposition}

\newtheorem{e-definition}[theorem]{Definition\rm}

\def\ch{{\mbox{ch}}}

\def\dim{{\mbox{dim}}}

\def\End{{\mbox{End}}}

\def\cala{{\mathcal A}}

\def\calr{{\mathcal R}}

\def\calw{{\mathcal W}}

\def\bbbone{\mbox{\rm 1\hspace {-.6em} l}}

\def\uw{\underline{W}}

\def\ug{{\mathbf u}}

\begin{document}

\thispagestyle{empty}
\begin{center}
{\bf GRADED ALGEBRAS AND MULTILINEAR FORMS}
\end{center} 

\begin{center}
 Michel DUBOIS-VIOLETTE
\footnote{Laboratoire de Physique Th\'eorique, UMR 8627, 
Universit\'e Paris XI,
B\^atiment 210\\ F-91 405 Orsay Cedex\\
Michel.Dubois-Violette$@$u-psud.fr}\\
\end{center}
\vspace{0,5cm}
 \vspace{0,5cm}

\begin{abstract}

We give a description of the connected graded algebras which are finitely generated and presented of global dimension 2 or 3 and which are Gorenstein.
These algebras are constructed from multilinear forms. We generalize the construction by associating homogeneous algebras to multilinear forms. The homogeneous algebras which are Koszul of finite global dimension and which are Gorenstein are of this type.
\end{abstract}

\section{Introduction}\label{sec1} One of our aims  is to study the connected graded algebras which are finitely generated in degree 1 and finitely presented with relations of degrees $\geq 2$ and which are of low global dimension $D$, $D=2$ and $D=3$. We further impose to these algebras to be Gorenstein. As pointed out in \cite{ber-mar:2006} (Proposition 5.2)  such an algebra is $N$-homogeneous and Koszul with $N=2$ for $D=2$ and $N\geq 2$ for $D=3$.  Our second more general objective is the study of the $N$-homogeneous algebras ($N\geq 2$) which are Koszul of arbitrary finite global dimension $D\geq 2$ and which are Gorenstein.  For $D=2$,  it is shown in  Section \ref{D2} that these algebras are classified by the nondegenerate bilinear forms modulo the action of the linear group.  For $D=3$, we show in Section \ref{D3} that these algebras are classified by an invariant subset of nondegenerate $(N+1)$-linear forms satisfying a regularity condition called 3-regularity modulo the action of the linear group.  In Section \ref{MF} we introduce and study homogeneous algebras associated with multilinear forms. It is pointed out that the Koszul homogeneous algebras of finite global dimension $D$ which are Gorenstein belong to this class which generalizes the previous results for $D=2,3$.\\
Throughout the paper $\mathbb K$ denotes a field which is algebraically closed, of characteristic zero and all the algebras and vector spaces are over $\mathbb K$. In the following $N, D$ and $q$ denote integers greater than or equal to 2 and we use the Einstein summation convention of repeated up down indices in the formulas. For the notion of Koszulity for homogeneous algebras introduced in \cite{ber:2001a} we refer to \cite{ber:2001a} and to \cite{ber-mdv-wam:2003}. For Koszul duality, Koszul $N$-complexes... for $N$-homogeneous algebras we refer to \cite{ber-mdv-wam:2003}; our notations are those of \cite{ber-mdv-wam:2003}.

\section{Global dimension $D=2$}\label{D2}
As explained in the introduction, the connected graded algebras which are finitely generated in degree 1 and finitely presented of global dimension 2 and which are Gorenstein are the quadratic Koszul algebras of global dimension 2 which are Gorenstein. These algebras are characterized by the following theorem.

\begin{theorem}\label{KG2}
Let $B$ be a nondegenerate bilinear form on $\mathbb K^q$ $(q\geq 2)$ with components $B_{\mu\nu}=B(e_\mu,e_\nu)$ in the canonical basis $(e_\lambda)_{\lambda\in \{1,\dots,q\}}$ of $\mathbb K^q$. Then the quadratic algebra $\cala$ generated by the elements $x^\lambda$ $(\lambda\in \{1,\dots, q\})$ with the relation $B_{\mu\nu} x^\mu x^\nu=0$ is Koszul of global dimension 2 and Gorenstein. Conversely, any quadratic algebra generated by $q$ elements $x^\lambda$ which is Koszul of global dimension 2 and Gorenstein is of the above kind for some non degenerate bilinear form $B$ on $\mathbb K^q$. Two nondegenerate bilinear forms on $\mathbb K^q$ which are on the same $GL(q,\mathbb K)$-orbit correspond to isomorphic algebras.
\end{theorem}

\noindent \underbar{Proof}. Let $\cala$ be the quadratic algebra generated by the $x^\lambda$ with the relation $B_{\mu\nu}x^\mu x^\nu=0$. Then the dual quadratic algebra $\cala^!$ is generated by elements $y_\lambda(\lambda\in \{1,\dots,q\})$ with relations $y_\mu y_\nu=(1/q) B_{\mu\nu} B^{\rho\tau}y_\tau y_\rho$ where $B^{\mu\nu}$ are the matrix elements of the inverse of the matrix $(B_{\mu\nu})$; i.e. $B_{\mu\lambda} B^{\lambda\nu}=\delta^\nu_\mu$. It follows that $\cala^!_0=\mathbb K\bbbone \simeq \mathbb K,\cala^!_1=\oplus_\lambda \mathbb K y_\lambda \simeq \mathbb K^q, \cala^!_2=\mathbb K B^{\beta\alpha} y_{\alpha} y_{\beta}\simeq \mathbb K$  and $\cala^!_n=0$ for $n\geq 3$. The Koszul complex of $\cala$ reads 
\begin{equation}
0\rightarrow \cala \stackrel{x^tB}{\rightarrow} \cala^q \stackrel{x}{\rightarrow} \cala \rightarrow 0
\label{KC2}
\end{equation}
where $\cala^q=(\cala,\dots, \cala)$, $x$ is right multiplication by the column $(x^\lambda)$ while $x^tB$ is right multiplication by the line $(x^\mu B_{\mu\nu})$. This complex is acyclic in degrees $\geq 1$ so the algebra $\cala$ is Koszul. The Gorenstein property follows by transposition and by using the invertibility of the matrix $(B_{\mu\nu})$.

Assume now that $\cala$ is a quadratic algebra generated by elements $x^\lambda$ $(\lambda\in \{ 1,\dots, q\})$ which is Koszul of global dimension 2 and Gorenstein. Then the Gorenstein property implies that the space of the quadratic relations of $\cala$ is of dimension 1, i.e. that the relations of $\cala$ read $B_{\mu\nu} x^\mu x^\nu=0$ for some non zero bilinear form $B$ on $\mathbb K^q$. The Koszul complex reads again as (\ref{KC2}) and the Gorenstein property implies that the matrix $(B_{\mu\nu})$ is invertible, i.e. that $B$ is nondegenerate. The fact that $\cala$ does  only depend on the $GL(q,\mathbb K)$-orbit of $B$ is straightforward.$\square$

For $q\geq 3$ the algebra has exponential growth while for $q=2$ it has polynomial growth so is regular in the sense of \cite{art-sch:1987}. In the latter case it is easy to classify the $GL(2,\mathbb K)$-orbits of nondegenerate bilinear forms on $\mathbb K^2$ according to the rank  of their symmetric part \cite{mdv-lau:1990} and one recovers the usual description of regular algebras of global dimension 2, \cite{irv:1979},  \cite{art-sch:1987}.

The algebra $\cala$ of Theorem \ref{KG2} corresponds to the natural quantum space for the action of the quantum group of the nondegenerate bilinear form defined in \cite{mdv-lau:1990}. 

\section{Global dimension $D=3$}\label{D3}
In order to state an analog of Theorem \ref{KG2} for homogeneous algebras which are Koszul of global dimension 3 and which are Gorenstein, we need the following concepts for multilinear forms.

\begin{e-definition}\label{PR}

Let $V$ be a vector space and $n$ be an integer with $n\geq 1$. A $(n+1)$-linear form $W$ on $V$ will be said to be preregular iff it satisfies the following conditions $(i)$, $(ii)$.

\noindent $(i)$ If $X\in V$ is such that one has  $W(X,X_1,\dots,X_n)=0$
for any $X_1,\dots, X_n \in V$, then $X=0$.

\noindent $(ii)$ There is an invertible linear transformation $Q_W\in GL(V)$ such that one has
\[
W(X_0,\dots,X_{n-1},X_n)=W(Q_W X_n,X_0,\dots, X_{n-1}) 
\]
for any $X_0,\dots,X_n \in V$.

\end{e-definition}

Let $W$ be a preregular $(n+1)$-linear form on $V$. Then the $Q_W\in GL(V)$ such that $(ii)$ is satisfied is unique and  given any $p\in \{0,\dots, n\}$ the condition $W(X_1,\dots,X_p,X,X_{p+1},\dots,X_n)=0$  for any $X_1,\dots,X_n \in V$
implies $X=0$. Furthermore the space of these $W$ is stable by the action $W\mapsto W^L$ of $GL(V)$ defined by
$W^L(X_0,\dots,X_n)=W(L^{-1}X_0,\dots,L^{-1}X_n)$ and one has
$Q_{W^L} = LQ_W L^{-1}$ for $L\in GL(V)$. Notice that a bilinear form is preregular iff it is nondegenerate.

\begin{e-definition}\label{3R}
Let $N$ be an integer with $N\geq 2$. A $(N+1)$-linear form $W$ on $V$ will be said to be 3-regular iff it is preregular and satisfies the following condition $(iii)$.\\
\noindent $(iii)$  If $L_0,L_1\in \End(V)$ are such that one has $W(L_0X_0,X_1,X_2,\dots,X_N)=W(X_0,L_1X_1,X_2,\dots,X_N)$ 
for any  $X_0,\dots,X_N \in V$, then $L_0=L_1=\lambda\bbbone$ for some $\lambda\in \mathbb K$.
\end{e-definition}
The subpace $\calr^N_3(V)$ of all 3-regular $(N+1)$-linear forms on $V$ is also stable by the action of $GL(V)$.

\begin{theorem}\label{KG3}
Let $\cala$ be a $N$-homogeneous algebra generated by $q$ elements $x^\lambda$ which is Koszul of global dimension 3 and Gorenstein. Then there is a 3-regular $(N+1)$-linear form $W$ on $\mathbb K^q$ with components $W_{\lambda_0\dots \lambda_N}$ in the canonical basis of $\mathbb K^q$ such that the relations of $\cala$ read $W_{\lambda\lambda_1\dots \lambda_N}x^{\lambda_1}{\cdots} x^{\lambda_N}=0$ ($\lambda\in\{1,\cdots,q\}$). Two 3-regular $(N+1)$-linear forms on $\mathbb K^q$ in the same $GL(q,\mathbb K)$-orbit correspond to isomorphic algebras.
\end{theorem}

\noindent {\bf Proof}.  
Given a $N$-homogeneous algebra $\cala$ which is Koszul of global dimension 3 and Gorenstein, the Koszul resolution of the trivial left $\cala$-module $\mathbb K$ reads
\[
0\rightarrow \cala\otimes (E\otimes R\cap R\otimes E)\stackrel{d}{\rightarrow} \cala\otimes R \stackrel{d^{N-1}}{\rightarrow} \cala\otimes E \stackrel{d}{\rightarrow} \cala\rightarrow \mathbb K\rightarrow 0
\]
and Gorenstein property implies  $\dim(E\otimes R\cap R\otimes E)=1$ so $E\otimes R\cap R\otimes E=\mathbb K W$ for some $W\in E^{\otimes^{N+1}}$ which (again via Gorenstein property) must satisfy conditions $(i)$ and $(ii)$,  $(iii)$ follows from exactness at $\cala\otimes R$. The fact that $\cala$ does only depend on the $GL(q,\mathbb K)$-orbit of $W$ is straightforward. $\square$\\

The Poincar\'e series of such a $N$-homogeneous algebra $\cala$ which is Koszul of global dimension 3 and which is Gorenstein is given by \cite{art-tat-vdb:1991},  \cite{mdv-pop:2002} $P_\cala(t)=(1-qt+qt^N-t^{N+1})^{-1}$
where $q=\dim (\cala_1)$  as before. It follows that $\cala$ has exponential growth if $q+N>5$. The case $q=2$ and $N=2$ is impossible so it remains the cases $q=3,N=2$ and $q=2,N=3$ for which one has polynomial growth \cite{art-sch:1987}. These latter cases are the object of \cite{art-sch:1987} and one encounters there various values of $Q_W$ in $GL(3,\mathbb K)$ and in $GL(2,\mathbb K)$. Examples with $N=3$ and arbitrary values of $q$ are the Yang-Mills algebra \cite{ac-mdv:2002b} and its deformations \cite{ac-mdv:2007} for which $Q_W=\bbbone$ and the super Yang-Mills algebra and its deformations \cite{ac-mdv:2007} for which $Q_W=-\bbbone$.

\section{Homogeneous algebras and multilinear forms}\label{MF}
Let $m$ and $N$ be two integers such that $m\geq N\geq 2$ and let $W$ be a $m$-linear form on $\mathbb K^q$ with $W\not=0$. We denote by $\uw$ the 1-dimensional space of multilinear forms spanned by $W$ and let $W_{\lambda_1\dots\lambda_m}=W(e_{\lambda_1},\dots,e_{\lambda_m})$ be the components of $W$ in the canonical basis $(e_\lambda)_{\lambda\in\{1,\dots, q\}}$ of $\mathbb K^q$. Let $\cala=\cala(W,N)$ be the $N$-homogeneous algebra generated by elements $x^\lambda$ with relations $W_{{\lambda_1}\dots{\lambda_{m-N}}{\mu_1}\dots {\mu_N}} x^{\mu_1}\dots x^{\mu_N}=0$ $(\lambda, \lambda_i\in \{1,\dots,q\})$,  $\cala(W,N)=A(E,R)$ with $E=\oplus_\lambda \mathbb K x^\lambda$ and   $R=\sum_{\lambda_i} \mathbb K W_{\lambda_1\dots \lambda_{m-N} \mu_1\dots \mu_N}x^{\mu_1} \otimes \dots \otimes x^{\mu_N}\subset E^{\otimes^N}$. The vector space $E$ can be interpreted as the dual of $\mathbb K^q$ and its basis $(x^\lambda)$ as the dual basis of $(e_\lambda)$ while $\uw$ is the 1-dimensional subspace of $E^{\otimes^m}$ spanned by $W\in E^{\otimes^m}$. To $\uw\subset E^{\otimes^m}$ we associate the family of ``derived" subspaces $\uw^{(n)}\subset E^{\otimes^{m-n}}$, $m\geq n\geq 0$, defined by $\uw^{(0)}=\uw$ and
$\uw^{(n)}=\sum_{\lambda_i} \mathbb K W_{\lambda_1\dots \lambda_n \mu_1\dots \mu_{m-n}} x^{\mu_1}\otimes \dots \otimes x^{\mu_{m-n}}$
for $m\geq n\geq 1$. Consider the sequence
\begin{equation}
0\rightarrow \calw_m \stackrel{d}{\rightarrow} \calw_{m-1}\stackrel{d}{\rightarrow} \dots \stackrel{d}{\rightarrow} \calw_2 \stackrel{d}{\rightarrow} \calw_1\stackrel{d}{\rightarrow}\calw_0\rightarrow 0
\label{wital}
\end{equation}
of free left $\cala$-modules and $\cala$-module-homomorphisms with $\calw_n\subset \cala\otimes E^{\otimes^n}$  defined by 
\[
\calw_n=\left \{ \begin{array}{lll}
\cala\otimes \uw^{(m-n)}& \mbox{if}\>  m\geq n\geq N\\
\cala\otimes E^{\otimes^n} & \mbox{if}\> N-1\geq n\geq 0
\end{array}\right.
\]
where the homomorphisms $d$ are induced by the homomorphisms of $\cala\otimes E^{\otimes^{n+1}}$ into $\cala\otimes E^{\otimes^n}$ defined by
$a\otimes (v_0\otimes v_1\otimes \dots v_n)\mapsto av_0\otimes (v_1 \otimes \dots \otimes v_n)$
for $m>n\geq 0$, $a\in \cala$ and $v_i\in E$ ($=\cala_1$).

\begin{e-proposition}\label{SubK}
Assume that the $m$-linear form $W$ is preregular. Then Sequence (\ref{wital}) is a sub-$N$-complex of the Koszul $N$-complex $K(\cala)$.
\end{e-proposition}

\noindent \underbar{Proof}. One has $W^{(m-n)}\subset E^{\otimes^{n-N}}\otimes R$ for $m\geq n\geq N$. By Property $(ii)$ of Definition \ref{PR}, relations 
$W_{\lambda_a\dots \lambda_{m-N}\mu_1\dots \mu_N} x^{\mu_1}\dots x^{\mu_N}=0$
are equivalent to relations
$W_{\lambda_{r+1}\dots \lambda_{m-N}\mu_1\dots \mu_N \nu_1\dots \nu_r}x^{\mu_1}\dots x^{\mu_N}=0$
for $m-N\geq r\geq 0$. It follows that $W^{(m-n)}\subset E^{\otimes^{n-N-r}}\otimes R\otimes E^{\otimes^r}$  so $\calw_n\subset \cala\otimes \cap_r E^{\otimes^{n-N-r}}\otimes R\otimes E^{\otimes^r}=K_n(\cala)$ for $m\geq n\geq N$ and one has the result.
$\square$

In the remaining part of this section $W$ is assumed to be preregular, $Q_W$ denotes the corresponding element of $GL(q,\mathbb K)=GL(E^\ast)$ and $\cala=\cala(W,N)=A(E,R)$ is the $N$-homogeneous algebra defined above. In view of the fact that one has $W\in \cap_r E^{m-N-r} \otimes R \otimes E^r=\cala^{!\ast}_m$, $W$ composed with the canonical projection $\cala^!\rightarrow \cala^!_m$ on degree $m$ defines a linear form $\omega_W$ on $\cala^!$.

\begin{e-proposition}\label{Aut}
Let $W,\cala,\cala^!,Q_W,\omega_W$ be as above.\\
$(i)$ $Q_W \in GL(E^\ast)$ induces an automorphism $\sigma_W$ of the $N$-homogeneous algebra $\cala^!$.\\
$(ii)$ One has $\omega_W(xy)=\omega_W(\sigma_W(y)x)$ for $x,y\in \cala^!$.
\end{e-proposition}

\noindent\underbar{Proof}. $Q_W$ induces an automorphism $\tilde \sigma_W$ of degree zero of the tensor algebra $T(E^\ast)$. Let $\tilde x$ be in $R^\perp$ that is
$\tilde x=\rho^{\mu_1\dots \mu_N} e_{\mu_1}\otimes \dots \otimes e_{\mu_N}$ is such that
$W_{\lambda_1\dots \lambda_{m-N} \mu_1\dots \mu_N} \rho^{\mu_1\dots \mu_N}=0, \>\> \forall \lambda_i$.
In view of $W\circ Q^{\otimes^m}_W=W$ it follows that one has
$Q^{\rho_1}_{\lambda_1}\dots Q^{\rho_{m-N}}_{\lambda_{m-N}}W_{\rho_1\dots \rho_{m-N} \nu_1\dots \nu_N}Q^{\nu_1}_{\mu_1}\dots Q^{\nu_N}_{\mu_N} \rho^{\mu_1\dots \mu_N}=0$ and therefore 
$W_{\lambda_1\dots \lambda_{m-N} \nu_1\dots \nu_N} Q^{\nu_1}_{\mu_1} \dots Q^{\nu_N}_{\mu_N} \rho^{\mu_1\dots \mu_N}=0$ ($ \forall \lambda_i$)
which means that $\tilde\sigma_W (\tilde x)$ is in $R^\perp$. Thus $\tilde\sigma_W$ passes to the quotient and defines an automorphism $\sigma_W$ of $\cala^!$ which proves $(i)$. 
$(ii)$ is then obvious in view of $(ii)$ in Definition \ref{PR}. $\square$

\begin{theorem}\label{KGD}
Let $\cala$ be a $N$-homogeneous Koszul algebra of finite global dimension $D$ which is Gorenstein. Then $\cala=\cala(W,N)$ for some preregular $m$-linear form $W$ on $\mathbb K^q$, $q=\dim (\cala_1)$ such that if $N\geq 3$ then $m=Np+1$ and $D=2p+1$ for some integer $p\geq 1$, while for $N=2$ one has $m=D$.
\end{theorem}
\noindent \underbar{Proof}. The Koszul resolution starts as $\dots \stackrel{d}{\rightarrow} \cala\otimes \cala^{!\ast}_N\stackrel{d^{N-1}}{\rightarrow} \cala\otimes\cala_1 \stackrel{d}{\rightarrow}\cala \rightarrow \mathbb K\rightarrow 0$ and must end as $0\rightarrow \cala \otimes \cala^{!\ast}_m\stackrel{d}{\rightarrow} \cala \otimes \cala^{!\ast}_{m-1}\stackrel{d^{N-1}}{\rightarrow}\dots$ in view of the Gorenstein property. This implies either $N=2$ and $m=D$ or if $N>2$, $m=Np+1$ and $D=2p+1 $ for some integer $p\geq 1$. The Gorenstein property implies also that $\dim(\cala^{!\ast}_m)=1$ and $\dim(\cala^{!\ast}_{m-1})=q$. Let $W$ be a generator of $\cala^{!\ast}_m\subset \cala^{\otimes^m}_1$. Identifying $\cala^\ast_1$ with $\mathbb K^q$, $W$ is a $m$-linear form on $\mathbb K^q$, $\cala^{!\ast}_{m-1}=\uw^{(1)}$ and $W$ satisfies Condition $(i)$ of Definition 3.1. Theorem 5.4 of \cite{ber-mar:2006} implies then that $\cala^{!\ast}_n=\uw^{(m-n)}$ for $n=Nk$ and $n=Nk+1$, $k\in \mathbb N$. So, in particular the space of relations $R=\cala^{!\ast}_N$ coincides with $\uw^{(m-N)}$ which means that $\cala=\cala(W,N)$. Condition $(ii)$ of Definition 3.1 follows also (see e.g. Corollary 5.12 of \cite{ber-mar:2006}), so $W$ is preregular.$\square$\\

Notice that, for even global dimension, a homogeneous Koszul-Gorenstein algebra is necessarily quadratic, (see in \cite{ber-mar:2006}, Proposition 5.3).
As example with $N=2$ and $m=D=4$, consider the algebra 
$\cala_\ug$ of \cite{ac-mdv:2002a} and \cite{ac-mdv:2003}. This is the quadratic algebra generated by 4 elements $x^\lambda$, $\lambda\in \{0,1,2,3\}$, with relations 
\[
\left\{
\begin{array}{l}
\cos(\varphi_0-\varphi_k)[x^0,x^k] = i\sin(\varphi_\ell-\varphi_m)\{x^\ell,x^m\} \\
\cos(\varphi_\ell-\varphi_m)[x^\ell,x^m] = i\sin(\varphi_0-\varphi_k) \{x^0,x^k\}
\end{array}
\right.
\]
for any cyclic permutation $(k,\ell,m)$ of (1,2,3) and where $\{A,B\}=AB+BA$. This algebra is Koszul of global dimension 4 and Gorenstein whenever none of these six relations becomes trivial and then one has as explained in \cite{ac-mdv:2002a}, the nontrivial Hochschild cycle 
\begin{eqnarray}
W  = \tilde{\ch}_{\frac{3}{2}}(U_\ug)= &- &\sum_{\alpha, \beta, \gamma,\delta} \epsilon_{\alpha \beta \gamma\delta} \cos(\varphi_\alpha-\varphi_\beta+\varphi_\gamma-\varphi_\delta)x^\alpha\otimes x^\beta\otimes x^\gamma\otimes x^\delta\nonumber \\
&+& i\sum_{\mu,\nu} \sin(2(\varphi_\mu-\varphi_\nu))x^\mu\otimes x^\nu\otimes x^\mu\otimes x^\nu\nonumber
\end{eqnarray}
which may be considered as a 4-linear form on $\mathbb K^4$. One checks that one has $\cala_\ug=\cala(W,2)$ and that $Q_W=-\bbbone$ here.

\end{document}